\documentclass{article}%
\usepackage{amsfonts}
\usepackage{graphicx}
\usepackage{amsmath}
\usepackage{amssymb}%
\providecommand{\U}[1]{\protect\rule{.1in}{.1in}}
\newtheorem{theorem}{Theorem}

\newtheorem{conjecture}[theorem]{Conjecture}
\newtheorem{corollary}[theorem]{Corollary}

\newtheorem{lemma}[theorem]{Lemma}

\newtheorem{problem}[theorem]{Problem}
\newtheorem{proposition}[theorem]{Proposition}

\newenvironment{proof}[1][Proof]{\textbf{#1.} }{\ \rule{0.5em}{0.5em}}
\begin{document}

\title{A much larger class of Fr\"{o}licher spaces than that of convenient vector
spaces may embed into the Cahiers topos}
\author{Hirokazu Nishimura\\Institute of Mathematics, University of Tsukuba\\Tsukuba, Ibaraki, 305-8571, Japan}
\maketitle

\begin{abstract}
It is well known that the category of Fr\"{o}licher spaces and smooth mappings
is Cartesian-closed. The principal objective in this paper is to show that the
full subcategory of Fr\"{o}licher spaces that believe in fantasy that every
Weil functor is really an exponentiation by the corresponding infinitesimal
object is also Cartesian-closed. Under the assumption that a conjecture holds,
it is shown that these Fr\"{o}licher spaces embed into the Cahiers topos with
the Cartesian-closed structure being preserved.

\end{abstract}

\section{Introduction}

It is widely held that the category of finite-dimensional smooth manifolds
does not behave itself well. It is not left exact, though it barely allows
transversal pullbacks. It is not cartesian closed, though it enjoys Weil
functors so that it admits infinitesimal objects to exponentiate over smooth
manifolds as the shadow of a shade. In order to develop really useful and
meaningful infinite-dimensional differential geometry, we should emancipate
ourselves from the regnant philosophy of \textit{manifolds}, which are
obtained by pasting together some parts of a linear space. The prototype of
the notion of manifold is traced back to Bernhart Riemann, and its modern
definition should be attributed to Hassler Whitney.

It is not difficult, after the manner of the well-established differential
geometry of finite-dimensional smooth manifolds, to introduce
infinite-dimensional manifolds by using Banach spaces, Hilbert spaces,
Fr\'{e}chet spaces or, more generally, convenient vector spaces in place of
finite-dimensional vector spaces. However the differential geometry of such
infinite-dimensional manifolds is neither so interesting theoretically nor so
useful practically. What is highly impressive, the last quarter of the 20th
century witnessed that Fr\"{o}licher, Kriegl and others finally arrived at the
correct notion of a smooth space, which is often called a Fr\"{o}licher space
in respect to his fame (cf. \cite{fro1}, \cite{fro2} and \cite{fro}). What is
most important and surprising about Fr\"{o}licher spaces is that they form a
cartesian closed category, while infinite-dimensional manifolds of any kind do not.

The notion of a \textit{Weil functor}, which is intended to stand for
exponentiation by the infinitesimal object corresponding to a Weil algebra in
fantasy, can be generalized easily to Fr\"{o}licher spaces in general. The
principal objective in this paper is to show that all the Fr\"{o}licher spaces
that believe that all Weil functors are really exponentiations by some
adequate infinitesimal objects in imagination form a cartesian closed
category. Such Fr\"{o}licher spaces are called \textit{Weil-exponentiable}. It
is then shown, under the assumption that a conjecture holds, that the
cartesian closed category of Weil-exponentiable Fr\"{o}licher spaces embeds
into the Cahier topos with the cartesian closed structure being preserved.
This is a direct generalization of Kock and Reyes' preceding result
(\cite{kock1} and \cite{kock2}) that the category of convenient vector spaces
embeds into the Cahier topos with the cartesian closed structure being
preserved. We hope that the paper will pave the way into a serious study on
the intriguing class of such modest Fr\"{o}licher spaces.

\section{Preliminaries}

\subsection{Fr\"{o}licher Spaces}

Fr\"{o}licher and his followers have vigorously and consistently developed a
general theory of smooth spaces, often called \textit{Fr\"{o}licher spaces}
for his celebrity, which were intended to be the central object of study in
infinite-dimensional differential geometry. A Fr\"{o}licher space is an
underlying set endowed with a class of real-valued functions on it (simply
called \textit{functions}) and a class of mappings from $\mathbb{R}$ to the
underlying set (called \textit{curves}) subject to the condition that curves
and functions should compose so as to yield smooth mappings from $\mathbb{R}$
to itself. It is required that the class of functions and that of curves
should determine each other so that each of the two classes is maximal with
respect to the other as far as they abide by the above condition. What is most
important among many nice properties about the category $\mathbf{FS}$ of
Fr\"{o}licher spaces and smooth mappings is that it is cartesian closed, while
the category of finite-dimensional smooth manifolds is not at all. It was
Lawvere et al. \cite{law} that established its cartesian closedness for the
first time, though their original proof has been simplified considerably. For
a standard reference on Fr\"{o}licher spaces the reader is referred to
\cite{fro}.

\subsection{Weil Algebras}

The notion of a \textit{Weil algebra} was introduced by Weil himself in
\cite{wei}. We denote by $\mathbf{W}$ the category of Weil algebras. Roughly
speaking, each Weil algebra corresponds to an infinitesimal object in the
shade. Although an infinitesimal object is undoubtedly imaginary in the real
world, each Weil algebra yields a \textit{Weil functor} on the category of
manifolds of some kind to itself. Intuitively speaking, the Weil functor
stands for the exponentiation by the infinitesimal object corresponding to the
Weil algebra at issue. For Weil functors on the category of finite-dimensional
smooth manifolds, the reader is referred to \S 35 of \cite{kolar}, while the
reader can find a readable treatment of Weil functors on the category of
smooth manifolds modeled on convenient vector spaces in \S 31 of \cite{kri}.

\subsection{Models of Synthetic Differential Geometry}

It is often forgotten that Newton, Leibniz, Euler and many other
mathematicians in the 17th and 18th centuries developed differential calculus
and analysis by using nilpotent infinitesimals. It is in the 19th century, in
the midst of the industrial revolution, that nilpotent infinitesimals were
overtaken by so-called $\varepsilon-\delta$ arguments. In the middle of the
20th century moribund nilpotent infinitesimals were revived by Grothendieck in
algebraic geometry and by Lawvere in differential geometry. Differential
geomery using nilpotent infinitesimals consistently is now called
\textit{synthetic differential geometry}, since it prefers synthetic arguments
to dull calculations as in ancient Euclidean geometry. Model theory of
synthetic differential geometry was developed vigorously by Dubuc and others
around 1980 by using techniques of topos theory. For a standard reference on
model theory of synthetic differential geometry the reader is referred to
Kock's \cite{kock}.

\section{Weil Prolongation}

The principal objective in this section is to assign, to each pair $(X,W)$\ of
a Fr\"{o}licher space $X$ and a Weil algebra $W$,\ another Fr\"{o}licher space
$X\otimes W$\ called the \textit{Weil prolongation of} $X$ \textit{with
respect to} $W$, which is naturally extended to a bifunctor $\mathbf{FS}%
\times\mathbf{W\rightarrow FS}$, and then to show that the functor
$\cdot\otimes W:\mathbf{FS\rightarrow FS}$ is product-preserving for any Weil
algebra $W$.

We will first define $X\otimes W$ set-theoretically for any pair $(X,W)$\ of a
Fr\"{o}licher space $X$ and a Weil algebra $W$. We define an equivalence
relation $\equiv\operatorname{mod}\ I$ on $\mathcal{C}^{\infty}(\mathbb{R}%
^{n},X)$\ to be
\begin{align*}
f  &  \equiv g\quad\operatorname{mod}\ I\\
&  \text{iff}\\
f(0,...,0)  &  =g(0,...,0)\text{ and}\\
\chi\circ f-\chi\circ g  &  \in I\text{ for any }\chi\in\mathcal{C}^{\infty
}(X,\mathbb{R})
\end{align*}
where $f,g\in\mathcal{C}^{\infty}(\mathbb{R}^{n},X)$. The totality of
equivalence classes with respect to $\equiv\operatorname{mod}\ I$ is denoted
by $X\otimes W$. This construction of $X\otimes W$ can naturally be extended
to a functor $\cdot\otimes W:\mathbf{FS\rightarrow Sets}$. We endow $X\otimes
W $ with the initial smooth structure with respect to the mappings
\[
X\otimes W\overset{\chi\otimes W}{\rightarrow}\mathbb{R}\otimes W
\]
where $\chi$ ranges over $\mathcal{C}^{\infty}(X,\mathbb{R})$. We note that
the smooth structure on $\mathbb{R}\otimes W$ has already been discussed by
Kock \cite{kock1}. Indeed Kock \cite{kock1} has indeed discussed how to endow
$X\otimes W$ with a smooth structure in case that $X$ is a convenient vector
space. Apparently the set-theoretical constructions of $\mathbb{R}\otimes W$
in this paper and Kock's \cite{kock1} coincide above all. Therefore we have to
verify that

\begin{proposition}
In case that $X$ is a convenient vector space, our above definition of the
smooth structure on $X\otimes W$ and that of Kock's \cite{kock1} coincide.
\end{proposition}

\begin{proof}
This follows readily from Theorem 1.3 of Kock's \cite{kock1}. The details can
safely be left to the reader.
\end{proof}

\begin{lemma}
\label{t2.1}Let $f,g\in\mathcal{C}^{\infty}(\mathbb{R}^{n},X)$ and $\varphi
\in\mathcal{C}^{\infty}(X,Y)$. Then
\[
f\equiv g\quad\operatorname{mod}\ I
\]
implies
\[
\varphi\circ f\equiv\varphi\circ g\quad\operatorname{mod}\ I
\]

\end{lemma}

Therefore each $\varphi\in\mathcal{C}^{\infty}(X,Y)$ defines a mapping
$X\otimes W\rightarrow Y\otimes W$ to be denoted by $\varphi\otimes W$, which
is to be shown smooth.

\begin{lemma}
\label{t2.2}The mapping $\varphi\otimes W:X\otimes W\rightarrow Y\otimes W$ is
really smooth.
\end{lemma}

\begin{proof}
This follows readily from the following commutative diagram:
\[%
\begin{array}
[c]{ccc}
& \mathbb{R}\otimes W & \\
\overset{(\chi\circ\varphi)\otimes W}{\nearrow} &  & \overset{\chi\otimes
W}{\nwarrow}\\
X\otimes W & \overset{\varphi\otimes W}{\rightarrow} & Y\otimes W\\
\uparrow &  & \uparrow\\
\mathcal{C}^{\infty}(\mathbb{R}^{n},X) & \overset{\mathcal{C}^{\infty
}(\mathbb{R}^{n},\varphi)}{\rightarrow} & \mathcal{C}^{\infty}(\mathbb{R}%
^{n},Y)
\end{array}
\]
where $\mathcal{C}^{\infty}(\mathbb{R}^{n},X)\rightarrow X\otimes W$ and
$\mathcal{C}^{\infty}(\mathbb{R}^{n},Y)\rightarrow Y\otimes W$ are the
canonical projections, and $\chi:Y\rightarrow\mathbb{R}$ is an arbitrary
smooth mapping.
\end{proof}

Let us assume that we are given two Weil algebras $W_{1}$ and $W_{2}$ of the
forms $\mathcal{C}^{\infty}(\mathbb{R}^{n})/I$ and $\mathcal{C}^{\infty
}(\mathbb{R}^{m})/J$ together with $\psi\in\hom_{\mathbf{W}}(W_{1},W_{2})$,
which can be represented by a smooth map $\overline{\psi}:\mathbb{R}%
^{m}\rightarrow\mathbb{R}^{n}\mathbb{\ }$abiding by the conditions:

\begin{enumerate}
\item $\overline{\psi}(\underset{m}{\underbrace{0,...,0}})=(\underset
{n}{\underbrace{0,...,0}})$

\item $\chi\in I$ implies $\chi\circ\overline{\psi}\in J$ for any $\chi
\in\mathcal{C}^{\infty}(\mathbb{R}^{n})$.
\end{enumerate}

Then $\overline{\psi}$ naturally induces a mapping $X\otimes W_{1}\rightarrow
X\otimes W_{2}$, which is denoted by $X\otimes\psi$.

\begin{lemma}
\label{t2.3}The mapping $X\otimes\psi:X\otimes W_{1}\rightarrow X\otimes
W_{2}$ is smooth.
\end{lemma}

\begin{proof}
This follows readily from the following commutative diagram:
\[%
\begin{array}
[c]{ccc}%
\mathbb{R}\otimes W_{1} & \overset{\mathbb{R}\otimes\psi}{\mathbb{\rightarrow
}} & \mathbb{R}\otimes W_{2}\\%
\begin{array}
[c]{cc}%
\chi\otimes W_{1} & \uparrow
\end{array}
&  &
\begin{array}
[c]{cc}%
\uparrow & \chi\otimes W_{2}%
\end{array}
\\
X\otimes W_{1} & \overset{X\otimes\psi}{\rightarrow} & X\otimes W_{2}\\
\uparrow &  & \uparrow\\
\mathcal{C}^{\infty}(\mathbb{R}^{n},X) & \overset{\mathcal{C}^{\infty
}(\overline{\mathbb{\psi}},X)}{\rightarrow} & \mathcal{C}^{\infty}%
(\mathbb{R}^{m},X)
\end{array}
\]
where $\mathcal{C}^{\infty}(\mathbb{R}^{n},X)\rightarrow X\otimes W_{1}$ and
$\mathcal{C}^{\infty}(\mathbb{R}^{m},X)\rightarrow X\otimes W_{2}$ are the
canonical projections, and $\chi:Y\rightarrow\mathbb{R}$ is an arbitrary
smooth mapping.
\end{proof}

By combining Lemmas \ref{t2.2} and \ref{t2.3}, we have

\begin{proposition}
\label{t2.4}We have the bifunctor $\otimes:\mathbf{FS}\times
\mathbf{W\rightarrow FS}$.
\end{proposition}

Now we are going to show that the functor $\cdot\otimes
W:\mathbf{FS\rightarrow FS}$ is product-preserving for any Weil algebra $W$.

\begin{lemma}
Let $\pi_{1}:X\times Y\rightarrow X$ and $\pi_{2}:X\times Y\rightarrow Y$ be
the canonical projections. The functions on $X\times Y$ are generated by the
compositions of $\pi_{1}$ followed by the functions on $X$ and the
compositions of $\pi_{2}$ followed by the functions on $Y$.
\end{lemma}

\begin{proof}
The reader is referred to Proposition 1.1.4 of \cite{fro}.
\end{proof}

\begin{corollary}
Let $W$ be of the form $\mathcal{C}^{\infty}(\mathbb{R}^{n})/I$. Let $X$ and
$Y$ be Fr\"{o}licher spaces. For any $f,g\in\mathcal{C}^{\infty}%
(\mathbb{R}^{n},X\times Y)$, we have that
\begin{align*}
f  &  \equiv g\quad\operatorname{mod}\ I\\
&  \text{iff}\\
\pi_{1}\circ f  &  \equiv\pi_{1}\circ g\quad\operatorname{mod}\ I\text{ }\\
&  \text{and }\\
\pi_{2}\circ f  &  \equiv\pi_{2}\circ g\quad\operatorname{mod}\ I
\end{align*}
where $\pi_{1}:X\times Y\rightarrow X$ and $\pi_{2}:X\times Y\rightarrow Y$
are the canonical projections.
\end{corollary}

\begin{proof}
It suffices to note that the product $X\times Y$ in the category $\mathbf{FS}
$ carries the initial structure induced by the mappings $\pi_{1}:X\times
Y\rightarrow X$ and $\pi_{2}:X\times Y\rightarrow Y$. More explicitly, the
functions on $X\times Y$ are generated by the compositions of $\pi_{1}$
followed by functions on $X$ and the compositions of $\pi_{2}$ followed by
functions on $Y$, for which the reader is referred to Proposition 1.1.4 of
\cite{fro}.
\end{proof}

\begin{corollary}
Let $U_{\mathbf{FS}}:\mathbf{FS\rightarrow Sets}$ be the forgetful functor.
Then the functor $U_{\mathbf{FS}}\circ(\cdot\otimes W):\mathbf{FS\rightarrow
Sets}$ is product-preserving, so that we have
\[
(X\times Y)\otimes W=(X\otimes W)\times(Y\otimes W)
\]
set-theoretically for any Fr\"{o}licher spaces $X$ and $Y$.
\end{corollary}

\begin{corollary}
The smooth structure of $(X\times Y)\otimes W$ is the initial structure with
respect to the canonical projections $(X\otimes W)\times(Y\otimes
W)\rightarrow X\otimes W$ and $(X\otimes W)\times(Y\otimes W)\rightarrow
Y\otimes W$.
\end{corollary}

Therefore we have

\begin{theorem}
Given a Weil algebra $W$, the functor $\cdot\otimes W:\mathbf{FS\rightarrow
FS}$ is product-preserving.
\end{theorem}

\section{Weil Exponentiability}

A Fr\"{o}licher space $X$ is called \textit{Weil-exponentiable }if
\begin{equation}
(X\otimes(W_{1}\otimes_{\infty}W_{2}))^{Y}=(X\otimes W_{1})^{Y}\otimes W_{2}
\label{4.1}%
\end{equation}
holds for any Fr\"{o}licher space $Y$ and any Weil algebras $W_{1}$ and
$W_{2}$. If $Y=1$, then (\ref{4.1}) degenerates into
\begin{equation}
X\otimes(W_{1}\otimes_{\infty}W_{2})=(X\otimes W_{1})\otimes W_{2} \label{4.2}%
\end{equation}
If $W_{1}=\mathbb{R}$, then (\ref{4.1}) degenerates into
\begin{equation}
(X\otimes W_{2})^{Y}=X^{Y}\otimes W_{2} \label{4.3}%
\end{equation}

\begin{proposition}
\label{t4.1}If $X$ is a Weil-exponentiable Fr\"{o}licher space, then so is
$X\otimes W$ for any Weil algebra $W$.
\end{proposition}

\begin{proof}
For any Fr\"{o}licher space $Y$ and any Weil algebras $W_{1}$ and $W_{2}$, we
have
\begin{align*}
&  ((X\otimes W)\otimes(W_{1}\otimes_{\infty}W_{2}))^{Y}\\
&  =(X\otimes((W\otimes_{\infty}W_{1})\otimes_{\infty}W_{2}))^{Y}\\
&  =(X\otimes(W\otimes_{\infty}W_{1}))^{Y}\otimes W_{2}\\
&  =((X\otimes W)\otimes W_{1})^{Y}\otimes W_{2}%
\end{align*}

\end{proof}

\begin{proposition}
\label{t4.2}If $X$ and $Y$ are Weil-exponentiable Fr\"{o}licher spaces, then
so is $X\times Y$.
\end{proposition}

\begin{proof}
For any Fr\"{o}licher space $Z$ and any Weil algebras $W_{1}$ and $W_{2}$, we
have
\begin{align*}
&  ((X\times Y)\otimes(W_{1}\otimes_{\infty}W_{2}))^{Z}\\
&  =\{\left(  X\otimes(W_{1}\otimes_{\infty}W_{2})\right)  \times\left(
Y\otimes(W_{1}\otimes_{\infty}W_{2})\right)  \}^{Z}\\
&  =\left(  X\otimes(W_{1}\otimes_{\infty}W_{2})\right)  ^{Z}\times\left(
Y\otimes(W_{1}\otimes_{\infty}W_{2})\right)  ^{Z}\\
&  =((X\otimes W_{1})^{Z}\otimes W_{2})\times((Y\otimes W_{1})^{Z}\otimes
W_{2})\\
&  =\{(X\otimes W_{1})^{Z}\times(Y\otimes W_{1})^{Z}\}\otimes W_{2}\\
&  =\{(X\otimes W_{1})\times(Y\otimes W_{1})\}^{Z}\otimes W_{2}\\
&  =\{(X\times Y)\otimes W_{1}\}^{Z}\otimes W_{2}%
\end{align*}

\end{proof}

\begin{proposition}
\label{t4.3}If $X$ is a Weil-exponentiable Fr\"{o}licher space, then so is
$X^{Y}$ for any Fr\"{o}licher space $Y$.
\end{proposition}

\begin{proof}
For any Fr\"{o}licher space $Z$ and any Weil algebras $W_{1}$ and $W_{2}$, we
have
\begin{align*}
&  (X^{Y}\otimes(W_{1}\otimes_{\infty}W_{2}))^{Z}\\
&  =(X\otimes(W_{1}\otimes_{\infty}W_{2}))^{Y\times Z}\\
&  =(X\otimes W_{1})^{Y\times Z}\otimes W_{2}\\
&  =((X\otimes W_{1})^{Y})^{Z}\otimes W_{2}\\
&  =(X^{Y}\otimes W_{1})^{Z}\otimes W_{2}%
\end{align*}

\end{proof}

\begin{theorem}
\label{t4.4}The full subcategory $\mathbf{FS}_{\mathbf{WE}}$\ of all
Weil-exponentiable Fr\"{o}licher spaces of the category $\mathbf{FS}$ of
Fr\"{o}licher spaces and smooth mappings is Cartesian closed.
\end{theorem}

\begin{problem}
Convenient vector spaces are all Weil-exponentiable. Find out a larger and
interesting class consisting of Weil-exponentiable Fr\"{o}licher spaces. If
any, find out a Fr\"{o}licher space which is not Weil-exponentiable.
\end{problem}

\section{The Embedding into the Cahiers Topos}

Let $\mathbf{D}$ be the full subcategory of the category of $\mathcal{C}%
^{\infty}$-rings of form $\mathcal{C}^{\infty}(\mathbb{R}^{n})\otimes_{\infty
}W$ with a natural number $n$ and a Weil algebra $W$. Now we would like to
extend the Weil prolongation $\mathbf{FS}_{\mathbf{WE}}\times\mathbf{W}%
\overset{\otimes}{\mathbf{\rightarrow}}\mathbf{FS}_{\mathbf{WE}}$ to a
bifunctor $\mathbf{FS}_{\mathbf{WE}}\times\mathbf{D}\overset{\otimes
}{\mathbf{\rightarrow}}\mathbf{FS}_{\mathbf{WE}}$. On objects we define
\begin{equation}
X\otimes C=\mathcal{C}^{\infty}(\mathbb{R}^{n},X)\otimes W \label{5.1}%
\end{equation}
for any Weil-exponentiable Fr\"{o}licher space $X$ and any $C=\mathcal{C}%
^{\infty}(\mathbb{R}^{n})\otimes_{\infty}W$. By Proposition \ref{t4.3}
$\mathcal{C}^{\infty}(\mathbb{R}^{n},X)$ is Weil-exponentiable, so that
$X\otimes C$ is Weil-exponentiable by Proposition \ref{t4.1}. It is easy to
see that the right hand of (\ref{5.1}) is functorial in $X$, but we have not
so far succeeded in establishing its functoriality in $C$. Therefore we pose
it as a conjecture.

\begin{conjecture}
\label{t5.1}The right hand of (\ref{5.1}) is functorial in $C$, so that we
have a bifunctor $\mathbf{FS}_{\mathbf{WE}}\times\mathbf{D}\overset{\otimes
}{\mathbf{\rightarrow}}\mathbf{FS}_{\mathbf{WE}}$.
\end{conjecture}

In the following we will assume that the conjecture is really true. We define
the functor $\mathbf{J:FS}_{\mathbf{WE}}\rightarrow\mathbf{Sets}^{\mathbf{D}}$
to be the exponential adjoint to the composite
\[
\mathbf{FS}_{\mathbf{WE}}\times\mathbf{D}\overset{\otimes}{\mathbf{\rightarrow
}}\mathbf{FS}_{\mathbf{WE}}\rightarrow\mathbf{Sets}%
\]
where $\mathbf{FS}_{\mathbf{WE}}\rightarrow\mathbf{Sets}$ is the
underlying-set functor. Now we have

\begin{proposition}
\label{t5.2}For any Weil-exponentiable Fr\"{o}licher space $X$ and any object
$\mathcal{C}^{\infty}(\mathbb{R}^{n})\otimes_{\infty}W$ in $\mathbf{D}$, we
have
\[
\mathbf{J}(X)^{\hom_{\mathbf{D}}(\mathcal{C}^{\infty}(\mathbb{R}^{n}%
)\otimes_{\infty}W,\cdot)}=\mathbf{J}(X\otimes(\mathcal{C}^{\infty}%
(\mathbb{R}^{n})\otimes_{\infty}W))
\]

\end{proposition}

\begin{proof}
Let $\mathcal{C}^{\infty}(\mathbb{R}^{m})\otimes W^{\prime}$ be an object in
$\mathbf{D}$. Then we have
\begin{align*}
&  \mathbf{J}(X)^{\hom_{\mathbf{D}}(\mathcal{C}^{\infty}(\mathbb{R}%
^{n})\otimes_{\infty}W,\cdot)}(\mathcal{C}^{\infty}(\mathbb{R}^{m}%
)\otimes_{\infty}W^{\prime})\\
&  =\hom_{\mathbf{sets}^{\mathbf{D}}}(\hom_{\mathbf{D}}(\mathcal{C}^{\infty
}(\mathbb{R}^{m})\otimes_{\infty}W^{\prime},\cdot),\mathbf{J}(X)^{\hom
_{\mathbf{D}}(\mathcal{C}^{\infty}(\mathbb{R}^{n})\otimes_{\infty}W,\cdot)})\\
&  =\hom_{\mathbf{sets}^{\mathbf{D}}}(\hom_{\mathbf{D}}(\mathcal{C}^{\infty
}(\mathbb{R}^{m})\otimes_{\infty}W^{\prime},\cdot)\times\hom_{\mathbf{D}%
}(\mathcal{C}^{\infty}(\mathbb{R}^{n})\otimes_{\infty}W,\cdot),\mathbf{J}%
(X))\\
&  =\hom_{\mathbf{sets}^{\mathbf{D}}}(\hom_{\mathbf{D}}(\mathcal{C}^{\infty
}(\mathbb{R}^{m+n})\otimes_{\infty}W\otimes_{\infty}W^{\prime},\cdot
),\mathbf{J}(X))\\
&  =\mathbf{J}(X)(\mathcal{C}^{\infty}(\mathbb{R}^{m+n})\otimes_{\infty
}W\otimes_{\infty}W^{\prime})\\
&  =\mathcal{C}^{\infty}(\mathbb{R}^{m+n},X)\otimes(W\otimes_{\infty}%
W^{\prime})\\
&  =\mathcal{C}^{\infty}(\mathbb{R}^{m},\mathcal{C}^{\infty}(\mathbb{R}%
^{n},X)\otimes W)\otimes W^{\prime}\\
&  =\mathbf{J}(X\otimes(\mathcal{C}^{\infty}(\mathbb{R}^{n})\otimes_{\infty
}W))(\mathcal{C}^{\infty}(\mathbb{R}^{m})\otimes_{\infty}W^{\prime})
\end{align*}

\end{proof}

\begin{proposition}
\label{t5.10}For any Weil-exponentiable Fr\"{o}licher spaces $X$ and $Y$, we
have
\[
\mathbf{J}(X\times Y)=\mathbf{J}(X)\times\mathbf{J}(Y)
\]

\end{proposition}

\begin{proof}
Let $\mathcal{C}^{\infty}(\mathbb{R}^{n})\otimes W$ be an object in
$\mathbf{D}$. Then we have
\begin{align*}
&  \mathbf{J}(X\times Y)(\mathcal{C}^{\infty}(\mathbb{R}^{n})\otimes_{\infty
}W)\\
&  =\mathcal{C}^{\infty}(\mathbb{R}^{n},X\times Y)\otimes W\\
&  =\{\mathcal{C}^{\infty}(\mathbb{R}^{n},X)\times\mathcal{C}^{\infty
}(\mathbb{R}^{n},Y)\}\otimes W\\
&  =\{\mathcal{C}^{\infty}(\mathbb{R}^{n},X)\otimes W\}\times\{\mathcal{C}%
^{\infty}(\mathbb{R}^{n},Y)\otimes W\}\\
&  =\mathbf{J}(X)(\mathcal{C}^{\infty}(\mathbb{R}^{n})\otimes_{\infty}%
W)\times\mathbf{J}(Y)(\mathcal{C}^{\infty}(\mathbb{R}^{n})\otimes_{\infty}W)\\
&  =(\mathbf{J}(X)\times\mathbf{J}(Y))(\mathcal{C}^{\infty}(\mathbb{R}%
^{n})\otimes_{\infty}W)
\end{align*}

\end{proof}

\begin{proposition}
\label{t5.3}For any Weil-exponentiable Fr\"{o}licher spaces $X$ and $Y$, we
have the following isomorphism in $\mathbf{Sets}^{\mathbf{D}}$:
\[
\mathbf{J}(X^{Y})=\mathbf{J}(X)^{\mathbf{J}(Y)}%
\]

\end{proposition}

\begin{proof}
Let $\mathcal{C}^{\infty}(\mathbb{R}^{n})\otimes W$ be an object in
$\mathbf{D}$. On the one hand, we have
\begin{align*}
&  \mathbf{J}(X^{Y})(\mathcal{C}^{\infty}(\mathbb{R}^{n})\otimes_{\infty}W)\\
&  =\mathcal{C}^{\infty}(\mathbb{R}^{n},X^{Y})\otimes W\\
&  =X^{Y\times\mathbb{R}^{n}}\otimes W\\
&  =(X^{\mathbb{R}^{n}}\otimes W)^{Y}%
\end{align*}
On the other hand, we have
\begin{align*}
&  \mathbf{J}(X)^{\mathbf{J}(Y)}(\mathcal{C}^{\infty}(\mathbb{R}^{n}%
)\otimes_{\infty}W)\\
&  =\hom_{\mathbf{Sets}^{\mathbf{D}}}(\hom_{\mathbf{D}}(\mathcal{C}^{\infty
}(\mathbb{R}^{n})\otimes_{\infty}W,\cdot),\mathbf{J}(X)^{\mathbf{J}(Y)})\\
&  =\hom_{\mathbf{Sets}^{\mathbf{D}}}(\hom_{\mathbf{D}}(\mathcal{C}^{\infty
}(\mathbb{R}^{n})\otimes_{\infty}W,\cdot)\times\mathbf{J}(Y),\mathbf{J}(X))\\
&  =\hom_{\mathbf{Sets}^{\mathbf{D}}}(\mathbf{J}(Y),\mathbf{J}(X)^{\hom
_{\mathbf{D}}(\mathcal{C}^{\infty}(\mathbb{R}^{n})\otimes_{\infty}W,\cdot)})\\
&  =\hom_{\mathbf{Sets}^{\mathbf{D}}}(\mathbf{J}(Y),\mathbf{J}(X\otimes
(\mathcal{C}^{\infty}(\mathbb{R}^{n})\otimes_{\infty}W)))\\
&  =\hom_{\mathbf{FS}_{\mathbf{WE}}}(Y,X\otimes(\mathcal{C}^{\infty
}(\mathbb{R}^{n})\otimes_{\infty}W))
\end{align*}

\end{proof}

\begin{theorem}
\label{t5.4}The functor $\mathbf{J}:\mathbf{FS}_{\mathbf{WE}}\rightarrow
\mathbf{Sets}^{\mathbf{D}}$ preserves the cartesian closed structure. In other
words, it preserves finite products and exponentials. It is full and faithful.
It sends the Weil prolongation to the exponentiation by the corresponding
infinitesimal object.
\end{theorem}

\begin{proof}
The first statement follows from Propositions \ref{t5.10} and \ref{t5.3}. The
second statement that it is full and faithful follows by the same token as in
\cite{kock1} and \cite{kock2}. The final statement follows from Proposition
\ref{t5.2}.
\end{proof}

The site of definition for the Cahiers topos $\mathcal{C}$ is the dual
category $\mathbf{D}^{\mathrm{op}}$\textbf{\ }of the category $\mathbf{D}$
together with the open-cover topology, so that we have the canonical
embedding
\[
\mathcal{C\hookrightarrow}\mathbf{Sets}^{\mathbf{D}}%
\]
By the same token as in \cite{kock1} and \cite{kock2} we can see that the
functor $\mathbf{J}:\mathbf{FS}_{\mathbf{WE}}\rightarrow\mathbf{Sets}%
^{\mathbf{D}}$ factors in the above embedding. The resulting functor is
denoted by $\mathbf{J}_{\mathcal{C}}$. Since the above embedding preserves
products and exponentials, Theorem \ref{t5.4} yields directly

\begin{theorem}
\label{t5.5}The functor $\mathbf{J}_{\mathcal{C}}:\mathbf{FS}_{\mathbf{WE}%
}\rightarrow\mathcal{C}$ preserves the cartesian closed structure. In other
words, it preserves finite products and exponentials. It is full and faithful.
It sends the Weil prolongation to the exponentiation by the corresponding
infinitesimal object.
\end{theorem}

\end{document}